\documentclass[latter,12pt, jhep]{article}
\usepackage{epsfig}
\usepackage{amssymb}
\usepackage{amsfonts}
\usepackage{amsmath}
\usepackage{mathrsfs}

\newskip\humongous \humongous=0pt plus 1000pt minus 1000pt

\newif\ifdtup

\catcode`\@=11

\@addtoreset{equation}{section}
\def\theequation{\thesection.\arabic{equation}}

\def\@normalsize{\@setsize\normalsize{15pt}\xiipt\@xiipt
\abovedisplayskip 14pt plus3pt minus3pt%
\belowdisplayskip \abovedisplayskip
\abovedisplayshortskip \z@ plus3pt%
\belowdisplayshortskip 7pt plus3.5pt minus0pt}

\def\small{\@setsize\small{13.6pt}\xipt\@xipt
\abovedisplayskip 13pt plus3pt minus3pt%
\belowdisplayskip \abovedisplayskip
\abovedisplayshortskip \z@ plus3pt%
\belowdisplayshortskip 7pt plus3.5pt minus0pt
\def\@listi{\parsep 4.5pt plus 2pt minus 1pt
     \itemsep \parsep
     \topsep 9pt plus 3pt minus 3pt}}

\relax

\catcode`@=12

\catcode`\@=11

\def\section{\@startsection{section}{1}{\z@}{3.5ex plus 1ex minus
   .2ex}{2.3ex plus .2ex}{\large\bf}}

\def\thesection{\arabic{section}}
\def\thesubsection{\arabic{section}.\arabic{subsection}}

\def\appendix{\setcounter{section}{0}
 \def\thesection{Appendix \Alph{section}}
 \def\thesubsection{\Alph{section}.\arabic{subsection}}
 \def\theequation{\Alph{section}.\arabic{equation}}}

\def\SymBoxes#1#2#3#4{\newdimen\un@t \un@t#3%
\raisebox{#1}{\rule{#2\un@t}{#4}\hskip-#2\un@t
\@tempdimb\un@t \advance\@tempdimb by-#4\@tempcntb#2\relax%
\@whilenum{\@tempcntb>0}\do{
\rule{#4}{\un@t}\hskip\@tempdimb \advance\@tempcntb by\m@ne}%
\hskip-#2\un@t \rule[\un@t]{#2\un@t}{#4}%
\rule[\un@t]{#4}{#4}\hskip-#4
\rule{#4}{\un@t}}\hskip-#4}                

\begin{document}

\newcommand{\beq}{\begin{equation}}
\newcommand{\eeq}{\end{equation}}
\newcommand{\bea}{\begin{eqnarray}}
\newcommand{\eea}{\end{eqnarray}}
\newcommand{\beas}{\begin{eqnarray*}}
\newcommand{\eeas}{\end{eqnarray*}}
\newcommand{\defi}{\stackrel{\rm def}{=}}
\newcommand{\non}{\nonumber}
\newcommand{\bquo}{\begin{quote}}
\newcommand{\enqu}{\end{quote}}
\renewcommand{\(}{\begin{equation}}
\renewcommand{\)}{\end{equation}}
\def\IZ{{\mathbb Z}}
\def\IR{{\mathbb R}}
\def\IC{{\mathbb C}}
\def\IQ{{\mathbb Q}}

\def\g{\gamma}
\def\m{\mu}
\def\n{\nu}
\def\a{\alpha}
\def\b{\beta}

\def\CM{{\mathcal{M}}}
\def\dCM{{\left \vert\mathcal{M}\right\vert}}

\def \d{\textrm{d}}
\def \p{\partial}

\def \Pf{\rm Pf\ }

\def \pr{\prime}

\def\Tr{ \hbox{\rm Tr}}
\def\half{\frac{1}{2}}

\def \eqn#1#2{\begin{equation}#2\label{#1}\end{equation}}
\def\de{\partial}
\def\Tr{ \hbox{\rm Tr}}
\def\H{ \hbox{\rm H}}
\def\HE{ \hbox{$\rm H^{even}$}}
\def\HO{ \hbox{$\rm H^{odd}$}}
\def\K{ \hbox{\rm K}}
\def\Im{ \hbox{\rm Im}}
\def\Ker{ \hbox{\rm Ker}}
\def\const{\hbox {\rm const.}}
\def\o{\over}
\def\im{\hbox{\rm Im}}
\def\re{\hbox{\rm Re}}
\def\bra{\langle}\def\ket{\rangle}
\def\Arg{\hbox {\rm Arg}}
\def\Re{\hbox {\rm Re}}
\def\Im{\hbox {\rm Im}}
\def\exo{\hbox {\rm exp}}
\def\diag{\hbox{\rm diag}}
\def\longvert{{\rule[-2mm]{0.1mm}{7mm}}\,}
\def\a{{\textsl a}}
\def\dag{{}^{\dagger}}
\def\tq{{\widetilde q}}
\def\p{{}^{\prime}}
\def\W{W}
\def\N{{\cal N}}
\def\hsp{,\hspace{.7cm}}
\newcommand{\C}{\ensuremath{\mathbb C}}
\newcommand{\Sp}{\ensuremath{\mathbb S}}
\newcommand{\Z}{\ensuremath{\mathbb Z}}
\newcommand{\R}{\ensuremath{\mathbb R}}
\newcommand{\rp}{\ensuremath{\mathbb {RP}}}
\newcommand{\cp}{\ensuremath{\mathbb {CP}}}
\newcommand{\vac}{\ensuremath{|0\rangle}}
\newcommand{\vact}{\ensuremath{|00\rangle}}
\newcommand{\oc}{\ensuremath{\overline{c}}}
\newcommand{\sgn}{\mathop{\mathrm{sgn}}}

\def\M{\mathcal{M}}
\def\F{\mathcal{F}}
\def\d{\textrm{d}}

\def\eps{\epsilon}

\begin{flushright}
\end{flushright}

\vspace{-.5truecm}
\begin{center}
{\Large \textbf{Convexity of a Small Ball Under Quadratic Map}}
\end{center}
\vspace{6pt}
\begin{center}
{\large\textsl{Anatoly Dymarsky $^{a,}$\footnote{Permanent address.}}\\}
\vspace{25pt}
\textit{ $^a$\small Center for Theoretical Physics, MIT, Cambridge, MA, USA, 02139}\\ \vspace{6pt}
\textit{ $^1$ Skolkovo Institute of Science and Technology,\\
Novaya St.~100, Skolkovo, Moscow Region, Russia, 143025} \vspace{6pt}
\end{center}

\vspace{6pt}
\begin{center}
\textbf{Abstract}
\end{center}
\vspace{-.1cm}
We derive an upper bound on the size of a ball such that the image of the ball under quadratic map is strongly convex and smooth. Our result is the best possible improvement of the analogous result by Polyak [1] in the case of quadratic map. We also generalize the notion of the joint numerical range of $m$-tuple of matrices by adding vector-dependent inhomogeneous terms and provide a sufficient condition for its convexity.
\vspace{6pt}
\begin{flushleft}
\end{flushleft}

\renewcommand{\thefootnote}{\arabic{footnote}}

\section{Introduction and Main Result}

\subsection{Polyak Convexity Principle}
\label{intro1}
Convexity is a highly appreciated feature which can drastically simplify analysis of various optimization and control problems. In most cases, however, the problem in question is not convex. In \cite{Polyak} Polyak proposed the following approach which proved to be useful in many applications \cite{PolyakApp}: to restrict the optimization or control problem to  a small convex subset of the original set.  More concretely, for a map  $y_i=f_i(x)$ from $\R^m$ to $\R^n$, instead of the full image ${\mathscr F}(f)\equiv f(\R^m)=\{f(x): x\in \R^n\}$, which is not necessarily convex, let us consider an image of a small ball $B_{\varepsilon}(x_0)=\{x: |x-x_0|^2\le \varepsilon^2\}$. For a regular point $x_0$ of $f_i(x)$ there is always  small $\varepsilon$ such that the image $f(B_{\varepsilon}(x_0))$ is convex. The underlying idea here is very simple: for any $x$ from a small vicinity of a  regular point $x_0$,  where ${\rm rank}\left({\partial f(x_0)/ \partial x}\right)=m$,  the map $f(x)$ can be approximated by a linear map
\bea
y_i(x)-y_i(x_0)\simeq \left. {\partial f_i\, \over \partial x^a}\right|_{x_0}(x-x_0)^a\ .
\eea
Since the linear map preserves strong convexity, so far the nonlinearities of $f(x)$ are small and can be neglected, the image of a small ball around $x_0$ will be convex. Reference \cite{Polyak} computes a conservative upper bound on $\varepsilon\le \varepsilon_P$ in terms of the smallest singular value $\nu$ of the Jacobian $J(x_0)\equiv \left. {\partial f \over \partial x}\right|_{x_0}$ and the Lipschitz constant $L$ of the Jacobian $\partial f(x)/\partial x$ inside $B_{\varepsilon}(x_0)$, 
\bea
\label{PR}
\varepsilon_P^2={\nu^2\over 4 L^2}\ .
\eea
The resulting image of $B_{\varepsilon}(x_0)$ satisfies the following two properties.\vspace{-.2cm}
\begin{enumerate}
  \item The image $f(B_{\varepsilon}(x_0))$ is strictly convex. 
  \item The pre-image of the boundary $\partial f(B_{\varepsilon}(x_0))$ belongs to the boundary $\partial B_{\varepsilon}(x_0)=\{x:   |x-x_0|^2=\varepsilon^2\}$. The interior points of $B_{\varepsilon}(x_0)$  are mapped into the interior points of $f(B_{\varepsilon}(x_0))$.
\end{enumerate}

\subsection{Local Convexity of Quadratic Maps}
In this paper we consider quadratic maps from $\R^n$ (or $\C^n$) to $\R^m$ of general form
\bea
\label{QM}
f_i(x)=x^* A_{i\,} x -v^*_i x-x^* v_i^{}\ ,
\eea
defined through an $m-$tuple of symmetric (hermitian) $n\times n$ matrices $A_i$ and an $m-$tuple of vectors $v_i\in \R^n$ (or $v_i\in \C^n$). Most of the results are equally applicable to both real $x\in \R^n$ and complex $x\in \C^n$ cases. The symbol ${}^*$ denotes transpose or hermitian conjugate correspondingly. Occasionally we will also use ${}^T$ to denote transpose for  the explicitly real-valued quantities. 

Applying the general theory of \cite{Polyak} toward \eqref{QM} one obtains 
\eqref{PR}, where $\nu^2$ is the smallest eigenvalues of the symmetric  $m\times m$ matrix $\Re(v^*_i v_j)$ and the Lipschitz constant $L$ for \eqref{QM} can be defined through
\bea
\label{Lipshitz}
L=\max_{|x^{}_1|^2=|x^{}_2|^2=1} \sqrt{ \sum\limits_{i=1}^n \Re(x_{1\,}^* A_{i\,} x^{}_2)^2}\ .
\eea
We see from \eqref{PR} that $\varepsilon_P$ is non-zero only if the point $x_0=0$ is regular and the $m\times 2 n$ matrix $\Re(v_i)\oplus \Im(v_i)$ has rank $m$.

Since any linear transformations of $x$ respects the form \eqref{QM}, generalizations to different central points $x_0\neq 0$ or  non-degenerate ellipsoids $(x-x_0)^*\, G\, (x-x_0)\le \varepsilon^2$, with some positive-definite $G$  instead of $|x-x_0|^2\le \varepsilon^2$ is trivial.

The bound \eqref{PR} is usually very conservative, and one can normally find a much larger ball $B_{\varepsilon'}(x_0)$ with $\varepsilon'>\varepsilon_P$ such that the properties $1, 2$ from section \ref{intro1} are still satisfied. The main result of this paper is the new improved bound $ \varepsilon^2_{\rm max}\ge \varepsilon^2_P$, where
\bea
\label{MR}
\varepsilon^2_{\rm max}\equiv\lim\limits_{\epsilon\rightarrow 0^+} \min\limits_{|c|^2=1} \left|(  c\cdot A-\lambda_{\rm m}(c\cdot A)+\epsilon)^{-1}c\cdot v\right|^2 \  ,\\
\lambda_{\rm m}(A)=\min\{\lambda_{\rm min}(A),0\}\ . \label{lm}
\eea
Here the minimum is over the unit sphere  from the dual space $c\in \R^m$,  $c\cdot y\equiv c^i y_i$, and $\lambda_{\rm min}(A)$ denotes the smallest eigenvalue of a symmetric (hermitian) matrix $A$. In \eqref{MR} and in what follows the sum of a matrix and a number always understood in a sense that the number is multiplied by $\mathbb I$, the  $n\times n$ identity matrix.

{\bf Proposition 1.} 
For any ball $B_\varepsilon (x_0=0)$, $\varepsilon^2< \varepsilon^2_{\rm max}$ the image $f(B_\varepsilon(0))$ is strongly  (strictly for $\varepsilon^2=\varepsilon^2_{\rm max})$ convex and smooth and the pre-image of the boundary $\partial f(B_\varepsilon(0))$ belongs to the boundary of the sphere $\partial B_\varepsilon(0)$ (properties $1, 2$ from section \ref{intro1}). The  value of \eqref{MR} is maximally possible  such that $f(B_\varepsilon(0))$ is stably convex (remains convex under infinitesimally small variation of $A_i, v_i$). In this sense \eqref{MR} is the best possible improvement of the Polyak's bound \eqref{PR}.

Sometimes propery $2$ is not important and can be relaxed. We would still want the image of $\partial B_\varepsilon(0)$ to be convex, but it is no longer important that the pre-image of the boundary $\partial f(B_\varepsilon(0))$ belongs solely to the boundary $\partial B_\varepsilon(0)$. In such a case the bound \eqref{MR} can be improved 
\bea
\label{MR2}
\tilde\varepsilon^2_{\rm max}\equiv\lim\limits_{\epsilon\rightarrow 0^+} \min\limits_{c\,\in\, {\mathcal C}} \left|(  c\cdot A-\lambda_{\rm min}(c\cdot A)+\epsilon)^{-1}c\cdot v\right|^2 \  ,\\
{\mathcal C}=\{c: c\in \R^m,\ |c|^2=1,\ \lambda_{\rm min}(c\cdot A)\le 0\}\ .
\eea 

{\bf Proposition 2.} 
For any ball $B_\varepsilon (x_0=0)$, $\varepsilon^2\le \tilde\varepsilon^2_{\rm max}$ the image $f(B_\varepsilon(0))$ is strictly convex (property $1$ from section \ref{intro1}).

To solve the minimization problems (\ref{MR}) and calculate the exact value of $\varepsilon_{\rm max}$ is a nontrivial task. In section \ref{CE} we introduce simplifications which lead to a number of easy-to-calculate conservative estimates of $\varepsilon_{\rm max}$. A reader primarily interested in practical applications of our results  can look directly there.  By further simplifying the bound \eqref{MR}  we recover \eqref{PR}. This is a first proof of the main result of \cite{Polyak} (in a particular case of quadratic maps) which is not based on the Newton's method. 

\section{Inhomogeneous Joint Numerical Range}
\label{IJNR}
Before discussing convexity of the a ball $B_\varepsilon(x_0)$ under quadratic map, it would be convenient first to understand the geometry of the image of a unit sphere $|x|^2=1$, 
\bea
\label{gjnr}
{\mathcal F}(A,v)=\left\{y_i: \exists\ x,\ y_i=f_i(x),\ |x|^2=1\right\}\ .
\eea  
Functions  $f_i(x)$ are defined in \eqref{QM}. We introduced a new notation   ${\mathcal F}(A,v)$ instead of the colloquial $f(|x|^2=1)$ to stress that  \eqref{gjnr} is an interesting object in its own right. We propose to call 
\eqref{gjnr} {\it inhomogeneous joint numerical range} because of its resemblance to the original definition: ${\mathcal F}(A,0)$ is the joint numerical range of the $m-$tuple of matrices $A_i$. Below we formulate a sufficient condition for ${\mathcal F}(A,v)$ to be convex. 

{\bf Proposition 3.} Inhomogeneous joint numerical range ${\mathcal F}(A,v)$ defined in \eqref{gjnr} is strongly convex and smooth if 
\bea
\label{suff}
\lim\limits_{\epsilon\rightarrow 0^+}\min_{|c|^2=1} \left| \left(c\cdot A -\lambda_{\rm min}(c\cdot A)+\epsilon\right)^{-1}c\cdot v\right| > 1\ ,
\eea
and $n>m$ ($2n>m$ in the complex case). When $n=2$ ($2n=m$) generalized joint numerical range ${\mathcal F}(A,v)$ is an ``empty shell''  ${\mathcal F}(A,v)=\partial_{\,}{\rm Conv}[{\mathcal F}(A,v)]$ where ${\rm Conv}$ denotes the convex hull.

{\it Comment}. The inequality \eqref{suff} is a sufficient  but not a necessary condition. Say, when all $v_i=0$,  inequality \eqref{suff} is not satisfied, but the joint numerical range ${\mathcal F}(A,0)$ can nevertheless be  (strongly) convex. This happens, for example, if the rank of the smallest eigenvalue of $c\cdot A$ is the same for all non-zero $c_i$ \cite{Gutkin} (see also \cite{Sheriff}). 

{\bf Proof of Proposition 3.}
We will prove strong convexity of ${\mathcal F}(A,v)$ using the supporting hyperplanes technique. Our logic closely follows the proof of convexity of the joint numerical range \cite{Gutkin}, \cite{Sheriff}.  Provided that \eqref{suff} is satisfied we will show that for any non-zero covector $c_i\in \R^m$ the corresponding supporting hyperplane touches ${\mathcal F}(A,v)$ at exactly one point.  This will establish  strict convexity of ${\rm Conv}[ {\mathcal F}(A,v)]$. By calculating the Hessian at the boundary and showing it is strictly positive  we establish that ${\rm Conv}[ {\mathcal F}(A,v)]$ is strongly convex.  Last, we  provide a topological argument that  for $n>m$ ($2n>m$ in the complex case) $f(x)$ is surjective inside $\partial {\mathcal F}(A,v)$.  

{\it Strict convexity of ${\rm Conv}[ {\mathcal F}(A,v)]$.}
Let us consider a non-zero covector $c_i \in \R^m$ and corresponding family of hyperplanes $c\cdot y ={\rm const}$. First, we would like to find a minimum $F_c=c\cdot y$ among all $y$ from ${\mathcal F}(A,v)$, which is the same as to minimize $c\cdot f(x)$ with the constraint $x^*x=1$. 
After introducing a Lagrange multiplier $\lambda$ to enforce the constraint, the equation determining $x$ takes the form 
\bea
\label{eqX}
(c\cdot A-\lambda)x=c\cdot v\ .
\eea 
Once $c_i$ is fixed it is convenient to diagonalize $c\cdot A$ and rewrite $c\cdot v$ using the eigenbasis 
\bea
(c\cdot A)x_k=\lambda_k\,  x_k,\quad x_k^*x^{}_l=\delta^{}_{kl}\ ,\\
c\cdot v=\sum\limits_{k=1}^n \alpha_k\, x_k\ .
\eea
Let us assume for now that all $\alpha_k\neq 0$.
The minimum of $F_c(x)$ is given by  a minimum of 
\bea
\label{fl}
F_c(\lambda)=\lambda-\sum\limits_k {|\alpha_k|^2\over \lambda_k-\lambda}\ ,
\eea 
subject to constraint 
\bea
\label{dF}
{d F_c\over d\lambda}=1-\sum\limits_k {|\alpha_k|^2\over (\lambda_k-\lambda)^2}=0\ .
\eea
In other words we need to find a local extremum of $F_c(\lambda)$ where $F_c(\lambda)$ is minimal.  This is not the same as the global minimum 
of $F_c(\lambda)$ because  this function is not bounded from below and approached minus infinity  when $\lambda \rightarrow -\infty$. In general the constraint ${d F_c/d\lambda}=0$ 
has many solutions (from $2$ to $2n$). Smallest $\lambda$ solving ${d F_c/d\lambda}=0$ corresponds to the smallest $F_c(\lambda)$.
Indeed, let $\tilde{\lambda}_1>\tilde{\lambda}_2$ be two solution of ${d F_c/d\lambda}=0$. Then 
\bea
F_c(\tilde{\lambda}_1)-F_c(\tilde{\lambda}_2)=\sum_{k} {|\alpha_k|^2 (\tilde{\lambda}_1-\tilde{\lambda}_2)^3\over 2(\lambda_k-\tilde\lambda_1)^2(\lambda_k-\tilde\lambda_2)^2}>0\ .
\eea The minimal $\lambda$ solving \eqref{dF} is smaller than all $\lambda_k$'s. On the interval from $-\infty$ to $\lambda_{\rm min}$ the derivative ${d F_c/d\lambda}$ monotonically decreases from $1$ to $-\infty$. Therefore it vanishes at exactly one point. At that point matrix $(c\cdot A -\lambda)$ is positive-definite and therefore $x$ minimizing  $c\cdot y(x)$ is unique.

Now we have to consider an  important case when  $\alpha_k=0$ for some $k=\tilde k$. Of course when $\lambda$ is generic $\lambda\neq \lambda_{\tilde k}$, function $F_c(\lambda)$ and the minimization problem remains the same. But in the special case $\lambda=\lambda_{\tilde k}$ matrix  $c\cdot A- \lambda$ develops a zero mode $x_{\tilde k}$ and the constraint $x^*x=1$ is no longer given by ${d F_c/d\lambda}=0$, but 
\bea
\label{dFc}
1\ge \sum\limits_{k\neq \tilde k} {|\alpha_k|^2\over (\lambda^{}_k-\lambda_{\tilde k})^2}=1-|x_{\tilde k\, }^* x|^2\ge 0\ .
\eea  
Comparing $F_c(\lambda)$ at two different solutions of \eqref{dFc} or \eqref{dF} we find that  $F_c$ is minimal at minimal $\lambda$. Hence if $c\cdot A$ has an eigenvalue $\lambda_{\tilde k}$ such that $x_{\tilde k\, }^* v=0$,  $\lambda_{\tilde k}$ satisfies \eqref{dFc}, and  it is smaller than the smallest solution of \eqref{dF} (which would imply $\lambda_{\tilde k}<\lambda^{}_k$ for all $k$ with $\alpha_k\neq 0$) the resulting  $x$ minimizing $c\cdot y(x)$ is not unique. This is because only the absolute value of the component of $x$ along $x_{\tilde k}$ is fixed by \eqref{dFc}, but not its sign (or the complex phase).   

Let us repeat what we understood so far.  For a given $c_i$, in case the projection of $c\cdot v$ on the eigenspace of the smallest eigenvalue of $c\cdot A$ is non-trivial, $(c\cdot v)^* x_{\rm min} \neq 0$, the supporting hyperplane orthogonal to $c_i$ always touches  ${\mathcal F}(A,v)$ ``from below'' at one point. In case $(c\cdot v)^* x_{\rm min}=0$ for all $x_{\rm min}$ corresponding to $\lambda_{\rm min}$, vector $x=(c\cdot A-\lambda_{\rm min}(c\cdot A)+\epsilon)^{-1}(c\cdot v)$ is well defined when $\epsilon\rightarrow 0^+$ and there are two options. If $|x|\ge 1$, the supporting hyperplane still touches  ${\mathcal F}(A,v)$ at one point, but if $|x|<1$ there are two (or continuum) points $x$ minimizing $F_c=c\cdot y(x)$ for $|x|^2=1$, although these point may still correspond to the same $y_i(x)$.

Above we explained that \eqref{suff} is a sufficient condition for 
${\rm Conv}[{\mathcal F}(A,v)]$ to be strictly convex. 

{\it The convex hull ${\rm Conv}[ {\mathcal F}(A,v)]$ is strongly convex and smooth. The boundary $\partial_{\,}{\rm Conv}[ {\mathcal F}(A,v)]$ is an embedding of $S^{m-1}$ in $\R^m$.} 
For $c_i\neq 0$
we define a map ${\rm y}_i(c)=y_i({\rm x}(c))$ where ${\rm x}(c)$ minimizes $c\cdot y(x)$ for $|x|^2=1$. As was demonstrated above such ${\rm x}(c)$ is unique for each $c$ and hence ${\rm y}(c)$ is well-defined. Since ${\rm y}(c)={\rm y}(\mu_{}c)$ for any positive $\mu$ function ${\rm y}(c)$ is defined on  the sphere $S^{m-1}\in \R^m$. First we will show that ${\rm y}:S^{m-1}\rightarrow \R^m$ is an immersion by proving that the rank of $\partial {\rm y}/\partial c$ is $m-1$ for all $c\neq 0$.  Let us introduce a ``time'' parameter $\tau$, such that $c_i(\tau=0)=c_i$ and $\dot{c}\equiv \left.{dc\over d\tau}\right|_{\tau=0}$ is a given vector from $TS^{m-1}$ which can be identified with  the orthogonal complement of $c$ inside $\R^m$. Given \eqref{suff} is satisfied there is a unique ${\rm x}(\tau)$ that minimizes $c(\tau)\cdot y(x)$ over $|x|^2=1$. At the point $\tau=0$ we have
\bea
\label{ma}
(c\cdot A-\lambda){\rm x}&=&c\cdot v\ , \\
(c\cdot A-\lambda)\dot{\rm x}&=&-(\dot{c}\cdot A-\dot\lambda){\rm x}+{\dot c}\cdot v\ .
\label{der}
\eea
Coefficient $\lambda(\tau)$ must be chose such  such that $|{\rm x}(\tau)|^2=1$, i.e.~${\rm x}$ is orthogonal to $\dot{\rm x}$. This is always possible because \eqref{der} becomes a non-degenerate linear equation on $\dot\lambda$ after multiplication by ${\rm x}^*(c\cdot A-\lambda)^{-1}$ from the left (as was discussed above matrix $(c\cdot A-\lambda)$ is positive-definite and hence non-degenerate). 

It follows from \eqref{der} that $\dot{\rm x}=0$ if and only if 
\bea
\label{ma2}
(\dot{c}\cdot A-\dot\lambda){\rm x}={\dot c}\cdot v\ .
\eea
We will show momentarily that \eqref{ma} combined with \eqref{ma2} would contradict the main assumption \eqref{suff}. 

{\bf Lemma 1.} If there is a vector $x$ of unit length, $|x|^2=1$, which solves 
\bea
(c_i\cdot A -\lambda_i)x=c_i\cdot v\quad {\rm for\quad }i=1,2\ ,
\eea
for two non-collinear $c_1, c_2$, $(c_1-c_2)\cdot v\neq0$ and two numbers $\lambda_1,\lambda_2$, and the matrix $c_1\cdot A-\lambda_1$ is positive-definite, then there exist $c\neq 0$, $\lambda$ such that $(c\cdot A -\lambda)$ is semi-positive definite, has a zero eigenvalue, and solves
\bea
(c\cdot A -\lambda)x=c\cdot v\ .
\eea
{\bf Proof of Lemma 1.}
Let us consider a one-dimensional family of vectors $c(\mu)=c_1(1+\mu)-{c}_2 \mu$ and function $\lambda(\mu)=\lambda_1(1+\mu)-\lambda_2 \mu$. Because $c_1, c_2$ are non-collinear  vector $c(\mu)\neq 0$ for any $\mu$. We know that 
\bea
\label{ma3}
(c(\mu)\cdot A-\lambda(\mu))x=c(\mu)\cdot v\ ,\quad |x|^2=1\ ,
\eea
for any $\mu$ and that for $\mu=0$ matrix $(c(\mu)\cdot A-\lambda(\mu))$ is positive-definite. When  $\mu\rightarrow \infty$, or $\mu\rightarrow -\infty$, or both,  matrix $(c(\mu)\cdot A-\lambda(\mu))$ will develop negative eigenvalues (unless $(c_1-c_2)\cdot A=\lambda_1-\lambda_2$; but then \eqref{ma3} would imply $(c_1-c_2)\cdot v=0$). Then by continuity there will be a value of $\mu$ when matrix  $(c(\mu)\cdot A-\lambda(\mu))$ is semi-positive definite with a zero eigenvalue. 

Now, using Lemma 1 with $c_1=c$ and $c_2=\dot{c}$ we find a contradiction with \eqref{suff}, which finishes our  proof that $\dot{\rm x}\neq 0$.  Hence $\rm{x}(c)$ is an immersion of $S^{m-1}$ into $S^{n-1}$ (or $S^{2n-1}$).

Finally we want to show that $\dot{\rm y} \neq 0$ for any $\dot{c}$ and hence ${\rm rank}(\partial {\rm y}/\partial c)=m-1$. It is enough to calculate 
\bea
\label{positive}
\dot{c} \cdot \dot{{\rm y}}=\dot{\rm x}^* (c\cdot A -\lambda(c)) \dot{\rm x}>0\ , 
\eea
where we used ${\rm x}^*{\rm \dot{x}}=0$ and positive-definiteness of $(c\cdot A -\lambda(c))$. Hence, ${\rm y}(c)$ is an immersion of $S^{m-1}$ into $\R^m$ and $\partial {\rm Conv}[{\mathcal F}(A,v)]$ is smooth. 

As a side note we observe that the second derivative of $c\cdot {\rm y}(\tau)$ is strictly positive as well
\bea
c\cdot \ddot{{\rm y}}=\dot{\rm x}^* (c\cdot A -\lambda(c)) \dot{\rm x}\ .
\eea Hence ${\rm Conv}[{\mathcal F}(A,v)]$ is strongly convex.  

To prove, that ${\rm y}(c)$ is an embedding we have to show that it is injective, i.e.~${\rm y}(c_1)={\rm y}(c_2)$ implies $c_1=c_2$. Clearly this would imply ${\rm x}={\rm x}(c_1)={\rm x}(c_2)$ as was discussed above. The vector ${\rm x}$ would satisfy 
\bea
(c_i\cdot A -\lambda_i)x=c_i\cdot v\quad {\rm for\quad }i=1,2\ ,
\eea
such that both matrices $(c_i\cdot A -\lambda_i)$ are positive-definite and therefore according to Lemma 1 this is inconsistent with \eqref{suff} unless $c_1=c_2$.

{\it Topological argument proving surjectivity of $y_i=f_i(x)$  on ${\rm Conv}[ {\mathcal F}(A,v)]$}.  Our last step is to show that ${\mathcal F}(A,v)$ coincides with its own convexification, i.e.~${\mathcal F}(A,v)$ includes all points contained inside $\partial_{\,} {\rm Conv}[ {\mathcal F}(A,v)]$. Let us assume this is not the case and there is a point $y_0$ in the interior of ${\rm Conv}[ {\mathcal F}(A,v)]$ which does not belong to ${\mathcal F}(A,v)$. Then  we can define a continuous retraction $\varphi$ of ${\mathcal F}(A,v)$ on $\partial_{\,}{\rm Conv}[ {\mathcal F}(A,v)]=S^{m-1}$. For any $y\in {\mathcal F}(A,v)$, we define $\varphi(y)$ as the intersection point of the ray from $y_0$ to $y$ and the  boundary $\partial_{\,}{\rm Conv}[ {\mathcal F}(A,v)]$. Because the set confined  by  the boundary $\partial_{\,}{\rm Conv}[ {\mathcal F}(A,v)]$ is convex $\varphi(y)$ is well-defined. Next, the embedding ${\rm y}({\rm x}(c))$ from $S^{m-1}$ to $\partial_{\,}{\rm Conv}[ {\mathcal F}(A,v)]$  can be inverted ${\rm y}^{-1}:\partial_{\,}{\rm Conv}[ {\mathcal F}(A,v)]\rightarrow S^{m-1}\subset S^{n-1}$ ($S^{m-1}\subset S^{2n-1}$ in the complex case). The combination $\phi={\rm y}^{-1}\circ \varphi \circ f$ defines a map $\phi: S^{n-1} \rightarrow S^{m-1}\subset S^{n-1}$ ( or $\phi: S^{2n-1} \rightarrow S^{m-1}\subset  S^{2n-1}$ in the complex case) from the sphere $|x|^2=1$ to the preimage of the boundary $\partial_{\,}{\rm Conv}[ {\mathcal F}(A,v)]$ inside $|x|^2=1$. Because $S^{m-1}$ is mapped by $\phi$ into itself it must be homologically non-trivial inside $S^{n-1}$ (or $S^{2n-1}$). This is possible only if $m=n$ ($m=2n$). In this case ${\mathcal F}(A,v)$ is an ``empty shell'', ${\mathcal F}(A,v)=\partial_{\,}{\rm Conv}[ {\mathcal F}(A,v)]$. Otherwise, when $n>m$ ($2n>m$), all $y_0$ confined by  $\partial_{\,}{\rm Conv}[ {\mathcal F}(A,v)]$ belong to ${\mathcal F}(A,v)$ which coincides with its convex hull. Because of \eqref{suff} matrix $v^a_i$ (or $\Re(v_i^a)\otimes\Im(v_i^a)$) must have rank $m$  which excludes $n<m$ ($2n<m$).

{\bf Corollary 2.} Sufficient condition \eqref{suff} can be relaxed to include equality. In such a case the set ${\mathcal F}(A,v)$ remains to be strictly convex when $n>m$ (or $2n>m$). If $n=m$  (or $2n=m$), ${\rm Conv}[{\mathcal F}(A,v)]$ remains to be strictly convex. 
At the same time \eqref{suff} can not be made much stronger. In case \eqref{suff} is larger than $1$, for some $c_i$ the minimum $F_c=c\cdot y(x)$ would be achieved at more than one $x$. Although the corresponding ${\mathcal F}(A,v)$ may still be convex, even strictly convex if all such $x$'s are mapped into the same point $y=f(x)$, convexity would be lost upon an infinitesimal variation $v_i\rightarrow v_i+\epsilon\, c^{\perp}_i x_{\rm min}$ (here $c^{\perp}_i$ is any vector orthogonal to $c_i$ and $\epsilon$ is an infinitesimal parameter).

\section{Convexity of Image of a Small Ball}
\subsection{Mapping  Interior Into Interior }
\label{III}
In this section we consider the main question, convexity of  the image of a ball $B_{\varepsilon}(0)=\{x: |x|^2\le \varepsilon^2\}$ under the map \eqref{QM} while the interior points of $B_{\varepsilon}(0)$ are mapped strictly into interior point of $f(B_{\varepsilon}(0))$. More precisely, we want to find a bound on $\varepsilon$ such that the image $f(B_{\varepsilon}(0))$ satisfies the properties $1, 2$ from section \ref{intro1}. 

It is convenient to think of the ball $B_{\varepsilon}(0)$ as a collection of spheres $|x|^2=z$, $\varepsilon^2\ge  z \ge0 $. Using results from the previous section we can easily find a point where each such sphere touches a supporting hyperplane defined by a vector $c\neq 0$. The minimal value of $c\cdot y(x)$  over $|x|^2=z$ is given by 
\bea
{\bf F}_{c}(z)=z\lambda-\sum_k {|\alpha_k|^2\over \lambda_k-\lambda}\ ,
\eea 
where $\lambda(z)$ is the smallest solution of (we assume  there is no $\lambda_{\tilde k}<\lambda(z)$, $\alpha_{\tilde k}=0$)
\bea
z=\sum_k {|\alpha_k|^2\over (\lambda_k-\lambda)^2}\ .
\eea 
From ${d{\bf F}_c\over dz}=\lambda$ we conclude that the ``outer layer'' of $f(B_{\varepsilon}(0))$ in the direction $c_i$ corresponds to $|x|^2=z$ with the maximal value of $z=\varepsilon^2$ if $\lambda(\varepsilon^2)<0$, or to $|x|^2=z<\varepsilon^2$ such that $\lambda(z)=0$ otherwise. In the latter case property 2 of section \ref{intro1} will not be satisfied: pre-image of $\partial f(B_{\varepsilon}(0))$ does not lie within the boundary $\partial B_{\varepsilon}(0)$. Combining the constraint that $f(|x|^2=\varepsilon^2)$ is the ``outer layer'' for all $c_i$ with the sufficient condition from section \ref{IJNR} ensuring its convexity  we obtain our main result $\varepsilon^2\le \varepsilon^2_{\rm max}$, where $\varepsilon^2_{\rm max}$ is defined in \eqref{MR} . 

For $n>m$ ($2n>m$) the topological argument from section \ref{IJNR} ensures that $f(|x|^2=\varepsilon^2)$ is a convex set and therefore $f(B_{\varepsilon}(0))=f(\partial B_{\varepsilon}(0))$. Hence Proposition 1 is proved in this case. For $n=m$ ($2n=m$) the image $f(|x|^2=\varepsilon^2)$ is an ``empty shell'' $f(\partial B_{\varepsilon}(0))=\partial f(B_{\varepsilon}(0))$ and we yet have to show that $f(B_{\varepsilon}(0))$ is convex by demonstrating that all points inside $\partial f(B_{\varepsilon}(0))$ belong to $f(B_{\varepsilon}(0))$. To that end we slightly modify the topological argument from section \ref{IJNR}. Let us assume there is $y_0$ inside $\partial f(B_{\varepsilon}(0))$ which does not belong to $f(B_{\varepsilon}(0))$. Then one can define a retraction $\varphi$ from  $f(B_{\varepsilon}(0))$  to $\partial f(B_{\varepsilon}(0))$. Hence we obtain the  map $\phi={\rm y}^{-1}\circ\varphi\circ f: B_{\varepsilon}(0) \rightarrow \partial f(B_{\varepsilon}(0))=S^{n-1}$, where possibility to invert ${\rm y}$  on $\partial f(B_{\varepsilon}(0))$ was proved in section  \ref{IJNR}. The boundary $S^{n-1}$ is mapped into itself by $\phi$ and therefore it is homologically non-trivial inside $B_{\varepsilon}(0)$, which is a contradiction. Hence such $y_0$ can not exist which finishes the proof of Proposition 1. 

\subsection{Mapping Interior Into Anything}
What if we relax property 2 of section \ref{intro1} and will no longer require the pre-image of $\partial f(B_{\varepsilon}(0))$ to belong solely to the boundary $\partial B_{\varepsilon}(0)$? We still would want to preserve strict convexity of $f(B_{\varepsilon}(0))$ (property 1). Recycling the results of section \ref{III} we conclude that  for any $c_i$ the corresponding supporting hyperplane  touches $f(B_{\varepsilon}(0))$ at exactly one point,  provided $\varepsilon^2\le \tilde\varepsilon_{\rm max}^2$, where $\tilde \varepsilon_{\rm max}$ is given by \eqref{MR2}.
Hence the ${\rm Conv}[ f(B_{\varepsilon}(0)]$ is strictly convex and $\partial_{\,}{\rm Conv}[ f(B_{\varepsilon}(0)]$ is  homeomorphic to a sphere $S^{m-1}$. 

Now the points from $\partial f(B_{\varepsilon}(0)$ may correspond not only to $\partial B_{\varepsilon}(0)$ but also to the interior of $B_{\varepsilon}(0)$. The remaining challenge is to modify the topological argument from the previous section to prove that all points confined by  $\partial_{\,}{\rm Conv}[ f(B_{\varepsilon}(0)]$ belong to $ f(B_{\varepsilon}(0))$.

Let $x={\rm x}(c)$ minimize $F_c\equiv c\cdot y(x)$ over $f(B_{\varepsilon}(0))$ for some $c_i\neq 0$. So far $\varepsilon \le \tilde \varepsilon_{\rm max}$ such ${\rm x}$ is unique for each ${\rm y}\in \partial_{\,}{\rm Conv}[ f(B_{\varepsilon}(0)]$. Hence we can define the continuous map ${\rm y}^{-1}$ from  $S^{m-1}=\partial_{\,}{\rm Conv}[ f(B_{\varepsilon}(0)]$ to $S^{m-1}\subset S^{n-1}$ (or $S^{m-1}\subset S^{2n-1}$).   
This is absolutely analogous to the cases considered previously, although now the map ${\rm y}^{-1}$ is merely a continuous homeomorphism, not an embedding as before. 

The rest is straightforward. Assuming there is $y_0$ inside $\partial_{\,}{\rm Conv}[ f(B_{\varepsilon}(0)]$ which does not belong to $ f(B_{\varepsilon}(0))$ we first form a retraction $\varphi: f(B_{\varepsilon}(0))\rightarrow \partial_{\,}{\rm Conv}[ f(B_{\varepsilon}(0)]$ and then the continuous map $\phi={\rm y}^{-1}\circ\varphi\circ f: B_{\varepsilon}(0) \rightarrow  S^{m-1}\subset B_{\varepsilon}(0)$. Since $S^{m-1}$ is mapped into itself it must be homologically non-trivial inside $B_{\varepsilon}(0)$ which is a contradiction. This finishes the proof of Proposition 2.

\section{Conservative Estimate of $\varepsilon^2_{\rm max}$}\label{CE}
Calculating $\varepsilon^2_{\rm max}$ from \eqref{MR} explicitly could be a difficult task. Therefore it would be useful to derive a conservative estimate $\varepsilon^2_{\rm est}\le \varepsilon^2_{\rm max}$ which would be easy to calculate in practice. For any given $c_i \neq 0$ the length of the vector  $( c\cdot A-\lambda_{\rm m}(c\cdot A)+\epsilon)^{-1}c\cdot v$ can be bound from below as follows
\bea
\label{est}
\lim\limits_{\epsilon\rightarrow 0^+} \left|(c\cdot A-\lambda_{\rm m}(c\cdot A)+\epsilon)^{-1}c\cdot v\right|^2\ge {|c\cdot v|^2\over || c\cdot A-\lambda_{\rm m}(c\cdot A)||^2 } \  .
\eea 
Here the matrix norm $|| A ||$ is defined as 
\bea
|| A ||\equiv \max_{|x|^2=1} |Ax|=\lambda_{\rm max}^{1/2}(A^* A)=\max\{\lambda_{\rm max}(A),\lambda_{\rm max}(- A)\}\ ,
\eea
where the last identity holds for a symmetric (hermitian) $A$.
In most cases the estimate \eqref{est} is very conservative. For example, if the projection of $c\cdot v$ on the zero eigenvector of $c\cdot A-\lambda_{\rm min}(c\cdot A)$ is non-vanishing, the LHS of \eqref{est} will be infinite while the RHS will stay finite. Nevertheless, the advantage of \eqref{est} is that it is much easier to deal with than  the original expression. 

The norm $|| c\cdot A-\lambda_{\rm m}(c\cdot A)||$ can be estimated  from above by $2|| c\cdot A|| $. This estimate is tight if $\lambda_{\rm max}(c\cdot A)=-\lambda_{\rm min}(c\cdot A)$ and is conservative otherwise. Hence we arrive at the following easy-to-calculate estimate 
\bea
\label{newb}
\varepsilon^2_{\rm est}=\min_{|c|^2=1}\, {|c\cdot v|^2\over 4 ||c\cdot A ||^2}\ .
\eea
This expression still can be simplified. To proceed further we would need the following lemma (below $||c_i||$ stands for the conventional definition of the vector norm  $||c_i||\equiv \sqrt{\sum\limits_{i=1}^n |c_i|^2}$).

{\bf Lemma 2.} For a $m$-tuple of symmetric (hermitian) matrices $A_i$
\bea
\label{lemma}
L(A)\equiv \max_{|x^{}_1|^2=|x^{}_2|^2=1} ||\Re( x_1^* A_{i\,} x^{}_2)||=\max_{
|x^{}_{}|^2=1} || x^* A_{i\,} x||=\max_{|c|^2=1} \lambda_{\rm max}(c\cdot A)\ .
\eea

{\bf Proof of Lemma 2.} For any symmetric (hermitian) matrix $A$, 
$\lambda_{\rm max}(A)=\max_{
|x^{}_{}|^2=1} (x^* A_{} x)$. Besides, for any real-valued vectors $a_i, c_i$, $\max_{|c|^2=1} (c\cdot a)=||a_i||$. Therefore 
\bea
\max_{|c|^2=1} \lambda_{\rm max}(c\cdot A)=\max_{|c|^2=1} \max_{
|x^{}_{}|^2=1} ( x^* (c\cdot A)_{} x)=\max_{
|x^{}_{}|^2=1} || x^* A_{i\, } x||\ ,
\eea
which proves the second equality of \eqref{lemma}.

It is obvious that $L(A)\equiv\max_{|x^{}_1|^2=|x^{}_2|^2=1} 
|| \Re(x_1^* A_{i\,} x^{}_2)||\ge \max_{
|x^{}_{}|^2=1} || x^* A_{i\,} x||$. Let $x_1^{\rm m}, x_2^{\rm m}$ be the vectors of unit length maximizing $||\Re(x_1^* A_{i\,} x^{}_2)||$. Let us also define a real-valued vector of unit length $c^{\rm m}_i=\Re((x_1^{\rm m})^* A_{i\,} x_2^{\rm m})/||\Re((x_1^{\rm m})^* A_{i\,} x_2^{\rm m})||$. Then 
\bea
L(A)=
\Re((x_1^{\rm m})^* (c^{\rm m}\cdot A)_{} x_2^{\rm m})\, \le\,  || c^{\rm m}\cdot A||\, \le\, \max_{|c|^2=1} \lambda_{\rm max}(c\cdot A)\ ,
\eea
which finished the proof. 

Let us now return back to \eqref{newb}. Because of $c_i$-dependence in both numerator and denominator this quantity may look difficult to calculate. Let us make one last  simplification and minimize/maximize numerator and denominator separately
\bea
\label{PR2}
\varepsilon^2_P={\min_{|c|^2=1} |c\cdot v|^2\over 4 \max_{|c|^2=1} \lambda^2_{\rm max}(c\cdot A)}\le \varepsilon^2_{\rm est}\ .
\eea
The obtained estimate $\varepsilon_P$ is nothing but the Polyak's bound (\ref{PR}, \ref{Lipshitz}). Hence we rederived Polyak's result in case of a  quadratic map without using Newton's method, something which has not been done before \cite{Polyak}. 

In fact we can do systematically better than \eqref{PR2}. Expression \eqref{newb} is homogeneous of zero degree in $c_i$ and therefore minimizing over $|c|^2=1$ or any other non-degenerate ellipse $c^*_{} g_{} c=1$ (where $g$ is a positive-definite $m\times m$ symmetric matrix) would yield the same result. 
This observation allows us to effectively get rid of the numerator of \eqref{newb} by means of preconditioning. We introduce matrix $g_{ij}=(v^*_i v_j+v^*_j v_i)/2$ and $\Lambda_i^j$ which transforms it into identity matrix, $\Lambda^T g\,\Lambda={\mathbb I}_{m\times m} $. Since $|c\cdot v|^2=c^*_{} g_{} c$ we obtain for \eqref{newb} 
\bea
\label{mymat}
\varepsilon^2_{\rm est}={1\over 4 L^2(\hat{A})}\ ,\qquad 
\hat{A}_i=\Lambda_i^j A_j\ .
\eea
 
\subsection{Approximate Estimate of $L$}
The problem of calculating or effectively estimating $L(A)$
is interesting in its own right. Originally Polyak provided an estimate 
\bea
L_P(A)=||\lambda_{\rm max}(A_i)||=\sqrt{\sum\limits_{i=1}^n ||A_i||^2}\ge L(A)\ ,
\eea
which is very conservative. Recently reference \cite{Xia} put forward a number of improvements, including a convex semidefinite programming algorithm which they claim accurately estimates $L(A)$ from above.  We believe the representation 
\bea
\label{lp}
L(A)^2=\max_{|x|^2=1}\, \sum\limits^n_{i=1}\, (x^* A_i^{}\,  x)^2\ ,
\eea
established in Lemma 2, would allow to reduce the problem of calculating $L(A)$ to one of the  known problems of convex optimization or provide the best possible effective algorithm to accurately estimate \eqref{lp}. Thus, by introducing the matrix $X=x\otimes x^*$ and relaxing the condition ${\rm rank}(X)=1$, \eqref{lp} can be recast as a minimization of a quadratic function over a convex space of positive-definite matrices. Furthermore, by treating $z=x\, \otimes\,  x$ as a vector in the $n^2$-dimensional space and introducing  $Z=z \otimes z^*$ after relaxing ${\rm rank}(Z)=1$ condition the problem reduces to a standard question of semi-definite programming. This trick was used in the algorithm of \cite{Xia}, which they conjectured to be the tightest estimate of $L(A)$ to date. We believe our method  will be more precise in estimating $L(A)$ because it is based on \eqref{lp}, which is an exact expression for $L(A)$, while the algorithm of \cite{Xia} relied on an approximate expression (also quartic in $x$) which bounds the true value of $L(A)$ from above. 

Besides the sophisticated algorithms to estimate $L(A)$ it would be of practical value to outline more elementary yet less precise ways to bound $\varepsilon^2_{\rm est}$. The original estimate by Polyak \eqref{lp} is very easy to calculate but it is too conservative. Reference \cite{Xia} suggests a systematic improvement over that result\footnote{The original paper \cite{Xia}   provides a slightly different formula $L_{\rm new} (A)\equiv \lambda_{\rm max}^{1/2}\left(\sum\limits_{i=1}^n A^*_i A^{}_i\right)$. This is of course the same for symmetric (hermitian) matrices $A_i$. A few examples  considered in \cite{Xia} include non-symmetric real matrices $A_i$. While one can define quadratic map \eqref{QM} with any real-valued $n\times n$ matrices $A_i$ (when $x\in \R^n$), only their symmetrization $(A_i^{}+A_i^T)/2$ truly  matters. Similarly the symmetrized matrices should be used to calculate $L_{\rm new}(A)$. Using non-symmetric $A_i$ would unnecessary increase the estimate $L_{\rm new}$.} 
\bea
\label{xia}
L_P(A)\,\ge\, L_{\rm new} (A)\equiv \lambda_{\rm max}^{1/2}\left(\sum\limits_{i=1}^n A^{2}_i\right)\,\ge\, L(A)\ .
\eea

Using last representation from Lemma 2 we can provide another estimate based on the  identity $\lambda^2_{\rm max}(c\cdot A)\le \Tr\left((c\cdot A)^2\right)$,
\bea
\label{myest}
L_{\rm n} (A)\equiv \lambda_{\rm max}^{1/2}\left(\Tr(A_{i\,} A_j)\right)\,\ge\, L(A)\ .
\eea
This bound can be further improved using identity (12) of \cite{Xia} (see there for original reference)
\bea
L_{\rm nov}&=&\max_{|c|^2=1}\, (c\cdot a +\sqrt{c^T M c})\ ,\\
a_i&=&\Tr(A_i)/m\ ,\qquad M_{ij}=\Tr(A_{i\,} A_j)-m\, a_i\, a_j\ .
\eea
Let $\lambda_1,\dots,\lambda_m$ be the eigenvalues of $M$ and ${\rm a}_k$ -- the projections of $a_i$ on the $k$-th eigenvector of $M$. Then calculating $L_{\rm nov}$ is reduced to finding roots of an algebraic equation 
\bea
&&L_{\rm nov}=\max_{\lambda\, \in\, {\mathcal L} } \sqrt{F_{\rm}(\lambda)}\ ,\qquad\qquad  \quad\, 
F_{\rm}(\lambda)=\lambda\left(1+\sum_k{ {\rm a}_k^2\over \lambda-\lambda_k}\right)\ ,\\
&&{\mathcal L}=\left\{\lambda: {d F_{\rm }(\lambda)\over d\lambda}=0\right\}\bigcup \left\{\lambda=\lambda_k: {\rm a}_k=0\ {\rm and}\ {d F_{\rm}(\lambda_k)\over d\lambda}>0\right\}\ .
\eea 
Neither $L_{\rm nov}$ nor $L_{\rm new}$ is systematically better. 

\section*{Acknowledgments}
I would like to thank Tudor Dimofte, Konstantin Turitsyn and Jamin Sheriff for useful discussions and gratefully acknowledge support from the grant RFBR 12-01-00482.


\end{document}